\def\misajour{07/03/2026}
\documentclass[english]{polyjmath}

\title{An introduction to the strategy \\ of transcendence proofs}
 
\author{Michel Waldschmidt}
\affil{Sorbonne Université, Université Paris Cité, CNRS, IMJ-PRG, F-75005 Paris, France\\
E-mail: { \tt \url{michel.waldschmidt@imj-prg.fr}}
\\
URL: {\tt \url{http://www.imj-prg.fr/~michel.waldschmidt}}}

\usepackage{amsmath,amsthm}

\usepackage{amssymb}
\usepackage{graphicx}
 \usepackage[applemac]{inputenc}
\usepackage{hyperref} 
\numberwithin{equation}{section}

\def\boxit#1#2{\setbox1=\hbox{\kern#1{#2}\kern#1}%
\dimen1=\ht1 \advance\dimen1 by #1 \dimen2=\dp1 \advance\dimen2 by #1
\setbox1=\hbox{\vrule height\dimen1 depth\dimen2\box1\vrule}%
\setbox1=\vbox{\hrule\box1\hrule}%
\advance\dimen1 by .4pt \ht1=\dimen1
\advance\dimen2 by .4pt \dp1=\dimen2 \box1\relax}
\def\hfl#1#2{\smash{\mathop{\hbox to 12 mm{\rightarrowfill}}
\limits^{\scriptstyle#1}_{\scriptstyle#2}}}

 \def\hlfl#1#2{\smash{\mathop{\hbox to 12 mm{\leftarrowfill}}
\limits^{\scriptstyle#1}_{\scriptstyle#2}}}
\def\phfl#1#2{\smash{\mathop{\hbox to 8 mm{\rightarrowfill}}
\limits^{\scriptstyle#1}_{\scriptstyle#2}}}

 \def\phlfl#1#2{\smash{\mathop{\hbox to 8 mm{\leftarrowfill}}
\limits^{\scriptstyle#1}_{\scriptstyle#2}}}
\def\atop#1#2{
\genfrac{}{}{0pt} {} 
{#1} 
{#2}}
\def\cqfd{\unskip\kern 6pt\penalty 500
\raise -2pt\hbox{\vrule\vbox to 10pt{\hrule width 4pt
\vfill\hrule}\vrule}\par}
\def\house#1{\setbox1=\hbox{$\,#1\,$}%
\dimen1=\ht1 \advance\dimen1 by 2pt \dimen2=\dp1 \advance\dimen2 by 2pt
\setbox1=\hbox{\vrule height\dimen1 depth\dimen2\box1\vrule}%
\setbox1=\vbox{\hrule\box1}%
\advance\dimen1 by .4pt \ht1=\dimen1
\advance\dimen2 by .4pt \dp1=\dimen2 \box1\relax}
\def\virgule{\raise 2pt \hbox{$,$}}
\def\C{\mathbb{C}} 
\def\Q{\mathbb{Q}}
 
\def\Z{\mathbb{Z}}

\def\og{\leavevmode\raise.3ex\hbox{$\scriptscriptstyle
\langle\!\langle$}}

\def\fg{\leavevmode\raise.3ex\hbox{$\scriptscriptstyle
\,\rangle\!\rangle$}}

 \newcommand\rmH{{\mathrm H}}
 \newcommand\rmd{{\mathrm d}}
 \newcommand\rme{{\mathrm e}}
 \newcommand\rmi{{\mathrm i}}
 
\begin{document}
 
\maketitle

\hfill {\em \small update: \misajour}
 
\begin{abstract}
The proof by Ch.~Hermite of the transcendence of the number $\rme$ in 1873 is the model on which most transcendence proofs for constants of analysis are based. The Schneider--Lang Theorem gives an opportunity to understand the strategy behind the proof. We only outline the main ideas, without providing technical details, for which we refer to the bibliography. This short paper does not contain any full proof. 
\end{abstract}

\medskip
\noindent
{\bf Accknowledgement} {\em
This note is a small part of a course on Transcendental Number Theory given in September 2024 at the AESIM School (Asian and European Schools In Mathematics), 
{\em Introduction to Number Theory and Algebraic Curves}, at the Mathematical Institute SASA, Faculty of Mathematics, University of Belgrade. The author is thankful to René Schoof who kindly suggested him to write down this part of the course. 
}

\section{Introduction}
Short proofs for the transcendence of numbers like $\rme$ or $\pi$ are available (see for instance \cite{Baker} for a one page proof for each of these two results, p.~4 for $\rme$, p.~5 for $\pi$). Our aim is not to deal with the challenge of producing proofs as concise as possible for these results as well as similar ones. Some readers may prefer to understand the ideas behind the arguments. In his groundbreaking paper \cite{Hermite}, Ch.~Hermite explains in a very pedagogical way what he is going to do. His very clever and original arguments rely on explicit formulae. 

In \cite{Siegel}, C.L.~Siegel obtained far reaching new results by avoiding explicit auxiliary results. This seminal contribution by Siegel paved the way to the solution by A.O.~Gel'fond \cite{Gelfond} and Th.~Schneider \cite{Schneider1} of Hilbert's 7th problem. 

A formalisation of Gel'fond's method was initiated by Th.~Schneider \cite{Schneider2} and considerably simplified by S.~Lang \cite{Lang}. The so--called {\em Schneider--Lang Theorem} (Theorem \ref{Theorem:SchneiderLang}) is a rather simple statement, based on Gel'fond's method, which implies many transcendence results including the theorems of Hermite--Lindemann and Gel'fond--Schneider. This Schneider--Lang Theorem offers a convenient way to explain the strategy of many transcendence proofs.

After Hermite's initial contribution to the subject \cite{Hermite}, very few papers dealing with transcendence proofs explained the basic ideas behind the technical arguments, until S.~Lang \cite{Lang}, not only contributed to the subject with new results dealing with algebraic groups (following a suggestion of P.~Cartier), but at the same time provided enlightening explanations on the proofs.

 We start (\S~\ref{S:HL}) with en explanation of what is going on in the proof by Gel'fond's method of the Hermite--Lindemann Theorem. Next (\S~\ref{S:G}) we show how a small modification of the same proof yields the Gel'fond--Schneider Theorem. Formalizing the ideas in \S~\ref{S:SL} following \cite[Ch.~III \S~1]{Lang}, we state the Schneider--Lang Theorem and sketch the proof.
 
A similar scheme of proof produces the solution of Hilbert's 7th problem via Schneider's method (\S~\ref{S:S}).

 Finally in \S~\ref{S:FurtherResults} we propose a short overview of other related results together with a list of references for going further in the subject.

\section{Hermite--Lindemann}\label{S:HL}
 The following result implies the transcendence of both numbers $\rme$ and of $\pi$ (using $\rme^{\rmi \pi}=-1$), as well as the transcendence of any nonzero logarithm of a nonzero algebraic number.

\begin{theorem}[Hermite--Lindemann]\label{Theorem:HL}
 If $\alpha$ is a nonzero complex number, then one at least of the two numbers $\alpha$, $\rme^{\alpha}$ is transcendental.
 \end{theorem}

This result concerns the exponential function $\rme^z$. There is no general statement on the transcendence of values of an analytic function, unless one requests special properties of this function. Let us give the list of specific properties of the exponential function 
\[
\rme^z=\sum_{n\geqslant 0} \frac {z^n}{n!} 
\]
 that we are going to need.

The first one is that it is an entire function, which means that it is analytic for all $z\in\C$. Also it is not a polynomial. For an entire function, it is equivalent to say that it is not a polynomial and to say that it is a transcendental function, namely a function which is not algebraic over the field of rational fractions in one variable: if $P$ is a nonzero polynomial in two variables, then the function $P(z,\rme^z)$ is not the zero function. Since the Taylor coefficients $1/n!$ of $\rme^z$ are rational numbers, it is equivalent to say that $\rme^z$ is transcendental over the field $\C(z)$ or over the field $\Q(z)$.
 
Next, the differential equation 
\[
\frac{\rmd}{\rmd z} \rme^z=\rme^z
\]
 plays an important role in the methods of Hermite and Gel'fond (not in Schneider's method, as we will see in \S\ref{S:S}).

Also the addition Theorem
\[
 \rme^{z_1+z_2}=\rme^{z_1}\rme^{z_2}
\]
 implies the multiplication formula
 \[
 \rme^{sz}=(\rme^z)^s
 \]
 for all $s\in\Z$, which will be used in all proofs below.

The last property that we will use is an upper bound for the growth of the entire functions $z$ and $\rme^z$: for $r\geqslant 0$,
\[
 \sup_{|z|=r}|z|=r
 \text{ and }
 \sup_{|z|=r}|\rme^z|=\rme^r.
 \]

Here is the strategy. We start with a nonzero complex number $\alpha$ and we assume that both $\alpha$ and $\rme^\alpha$ are algebraic. Our goal is to get a contradiction with the fact that the exponential function is transcendental: we wish to prove that there exists a nonzero polynomial $P\in\Z[X,Y]$ such that the function $F(z)=P(z,\rme^z)$ is the zero function. We will achieve this goal in several steps. The idea is to show the existence of such a polynomial for which the associated function $F$ has many zeroes. We will start by requesting only a finite number of zeroes, so that linear algebra suffices to ensure the existence of $P$. A refinement of the argument of linear algebra is needed, because we will need an upper bound not only on the degrees of the polynomial $P$ but also on its coefficients. This refinement rests on Dirichlet's box principle, by means of the so--called Thue--Siegel Lemma.

 In order to achieve this first step, one will need to introduce some parameters which will be upper bounds for the degrees in $X$ and $Y$ of the polynomial $P$, and also on the number of zeroes that we are requesting.

The information we have is that $\alpha$ and $\rme^\alpha$ are algebraic. For a polynomial $P$ in two variables with integer coefficients, the values of the function $F(z)=P(z,\rme^z)$ at the point $z=\alpha$ is algebraic, in the number field $K:=\Q(\alpha,\rme^\alpha)$. Thanks to the differential equation, also the derivatives of $F$ at this point are also in $K$. Further, thanks to the multiplication formula, all numbers
 \begin{equation}\label{Equation:Fsigma(salpha)}
 F^{(\sigma)} (s\alpha) :=
\left(\frac{\rmd}{\rmd z} \right)^\sigma F(s\alpha) 
\end{equation}
for $\sigma\geqslant 0$ and $s\in\Z$ are in $K$. These are the numbers we are interested in. The first step will be requesting that a finite set of such numbers \eqref{Equation:Fsigma(salpha)} vanish. 

After this first step, we will increase the set of $(\sigma,s)$ for which these numbers \eqref{Equation:Fsigma(salpha)} vanish. This second step relies on the fact that a function having many zeroes is {\em flat}, it takes small values in a disc containing the given zeroes, and the same is true for the first derivatives. This analytic argument based on Schwarz's Lemma will provide an upper bound for the absolute value of numbers of the form \eqref{Equation:Fsigma(salpha)}. 

Now comes the arithmetic estimate, one of Liouville's inequalities. The idea is that a nonzero rational integer has an absolute value at least $1$, a nonzero rational number with denominator $d$ has an absolute value at least $1/d$, and more generally one can give a lower bound for the absolute value of a nonzero algebraic number when one knows an upper bound for its conjugates and denominator 
 (see \cite[Ch.~III]{LILAN}, \cite[\S~1.1]{Queens}, \cite[Ch.~3]{GL326}). It follows that these algebraic numbers \eqref{Equation:Fsigma(salpha)} all vanish. The contradiction follows at once. 

Here is a slightly expanded description of the proof. 

One needs to estimate certain quantities (something we only outline here). This will rely on introducing auxiliary parameters which will be supposed to be very large, so that only asymptotic estimates will suffice. 

\medskip
\noindent
{\bf Step 0: \em introducing parameters}
\\
We will show the existence of an auxiliary polynomial $P\in \Z[X,Y]$ satisfying suitable properties; we introduce upper bounds for its degrees, namely $T_0$ for the degree in $X$ and $T_1$ for the degree in $Y$. Hence the polynomial $P$ will be of the form
\[
P(X,Y)=\sum_{\tau=0}^{T_0-1}\sum_{t=0}^{T_1-1}a_{\tau,t}p_{\tau,t}(X,Y)
\]
with $a_{\tau,t}\in\Z$, where $p_{\tau,t}(X,Y)$ is the monomial $X^\tau Y^t$.
The maximum of the absolute value of its coefficients will be denoted by $\rmH(P)$. We will consider the entire function 
\[
F(z)=P(z,\rme^z)=\sum_{\tau=0}^{T_0-1}\sum_{t=0}^{T_1-1} a_{\tau,t}f_{\tau,t}(z)
\text{ 
where
}
f_{\tau,t}(z):=p_{\tau,t}(z,\rme^z)=z^\tau \rme^{tz}.
\] 
Its derivatives at points multiple of $\alpha$ are
\[
F^{(\sigma)} (s\alpha)=
\sum_{\tau=0}^{T_0-1}\sum_{t=0}^{T_1-1} a_{\tau,t}
f^{(\sigma)}_{\tau,t}(s\alpha)
\]
with $\sigma\geqslant 0$ and $s\geqslant 0$.

An explicit formula for the derivatives of $f_{\tau,t}$ is 
\[
f^{(\sigma)} _{\tau,t}(z)=
\sum_{\kappa=0}^{\min\{\tau,\sigma\}}
\frac
{\sigma!}
{\kappa!(\sigma-\kappa)!}
\cdot
\frac
{\tau!}
{(\tau-\kappa)!} t^{\sigma-\kappa}
z^{\tau-\kappa} \rme^{tz},
\]
which shows that the numbers 
\begin{equation}\label{Equation:Coefficients}
f^{(\sigma)} _{\tau,t}(s\alpha)=
\sum_{\kappa=0}^{\min\{\tau,\sigma\}}\frac
{\sigma!\tau!}
{\kappa!(\sigma-\kappa)!(\tau-\kappa)!}
t^{\sigma-\kappa}
(s\alpha)^{\tau-\kappa}(\rme^\alpha)^{st}
\end{equation}
for $\sigma \geqslant 0$, $s\geqslant 0$, $\tau\geqslant 0$, $t\geqslant 0$ all belong to $K$.

When writing down the proof, one first select parameters $S_0,S_1,T_0,T_1$ so that the proof can be carried out. Here we will introduce the contraints on these parameters only when they arise. One could produce similar sketches of proofs for many open problems, only the constraints which will occur during the process will turn out at the end not to be compatible. 

We denote by $\delta$ the degree of the number field $K=\Q(\alpha,\rme^\alpha)$. This number $\delta$ is a constant, while $S_0$, $S_1$, $T_0$, $T_1$ can be arbitrarily large. 

\medskip
\noindent
{\bf Step 1: \em existence of an auxiliary function}
\\
With this method one does not {\em construct } an auxiliary function, we merely prove that it exists, based on Dirichlet's box principle. 

Assume first $T_0T_1>S_0S_1$. Then linear algebra implies the existence of numbers $ a_{\sigma,s}$ in $K$, ($0\leqslant \sigma<S_0$, $0\leqslant s<S_1$), not all of which are $0$, such that 
\begin{equation}\label{Equation:FirstStep}
F^{(\sigma)} (s\alpha) =0
\quad\text{ for }\quad 0\leqslant \sigma<S_0,\; 0\leqslant s<S_1.
\end{equation}
If we assume the stronger condition $T_0T_1>\delta S_0S_1$, then there is a solution to this system of homogeneous linear equations with $ a_{\sigma,s}$ in $\Z$. This is not quite sufficient: we will need an estimate for the height 
\[
\rmH(P)=\max_{\atop{0\leqslant \tau<T_0}{0\leqslant t<T_1}}|a_{\tau,t}|
\]
of $P$. Estimates for the solutions of the linear system given by linear algebra would not suffice. One main idea introduced by A.~Thue and developed by C.L.~Siegel \cite{Siegel} is that a sharp upper bound can be achieved if one requests
\begin{equation}\label{Equation:T0T1S0S1}
T_0T_1\geqslant 2 \delta S_0S_1.
\end{equation}
This is our first constraint on the parameters. 

The coefficients of the system of linear equations are given by \eqref{Equation:Coefficients}. They belong to $K$. The version of the Thue-Siegel Lemma in \cite{Siegel} deals with the case of rational integer coefficients. One easily reduces to this case by writing the powers of $\alpha$ and $\rme^\alpha$ on a basis of $K$ over $\Q$ and multiplying by a denominator. Variants of the Thue--Siegel Lemma are given in \cite[\S~1.3]{LN402} and \cite[\S~1.2]{Queens}. 

\medskip
\noindent
{\bf Step 2: \em second step of the induction}
\\
Step 1 is the first step of the inductive argument. We explain the second step before performing the induction. We wish to increase the number of equations \eqref{Equation:FirstStep}. We could increase the range of $s$, but we choose to increase the range or $\sigma$. Hence we fix $s$ satisfying $0\leqslant s<S_1$ and we consider the number 
\[
\gamma:= F^{(S_0)} (s\alpha).
\]
This number belongs to $K$. By Cauchy's estimate, it satisfies 
\[
|\gamma|\leqslant S_0! |F|_r
\]
for any $r\geqslant 1+s|\alpha|$, where
\[
|F|_r:=\sup_{|z|=r}|F(z)|.
\]
Now the entire function $F$ has at least $S_0S_1$ zeroes (counting multiplicities) in the disc $|z|\leqslant r$. Here comes the analytic estimate, which is Schwarz's Lemma (\cite[\S~1.3]{Queens}, \cite[Lemma 7.1.3]{Asterisque}): an analytic function having many zeroes, say $N$, in a disc $|z|\leqslant r$, takes only small values in this disc: 
\[
|F|_r\leqslant \left( \frac R {3r}\right)^{-N}|F|_R
\]
for $R>3r$. The factor $(R/3r)^{-N}$ is small when $N$ (here $N\geqslant S_0S_1$) is large and $R>3r$. Therefore $|\gamma|$ is small. 

Now comes the arithmetic estimate, namely the above mentioned Liouville's inequality: the upper bound for the algebraic number $|\gamma|$ implies $\gamma=0$. 

Once we are able to increase the number of zeroes of our auxiliary function $F$, we are in a position of performing an induction which at the end will imply the desired conclusion: $F$ is the zero function. 

Our first constraint on the parameters $T_0,T_1, S_0, S_1$ occurred in step 1. The estimates necessary to complete the proof of step 2 introduce further constraints which lead to the choice of these parameters. A further quantity has been introduced, namely the large radius $R$, which needs also to be selected in terms of the parameters $T_0,T_1, S_0, S_1$.

See \cite[Ch.~3]{LN402} for all details, including a suitable choice of the parameters.

\medskip
\noindent
{\bf Step 3: \em the inductive argument}
\\
There is no logical need to perform the second step separately, one can directly perform the induction and prove that for any $\sigma\geqslant 0$ and for $1\leqslant s<S_1$ we have
\[
F^{(\sigma)} (s\alpha)=0.
\]
For $0\leqslant \sigma<S_0$, these equations are proved in the first step. For $\sigma=S_0$, this is the second step above. Assume this is true for all $\sigma<S$ for some $S\geqslant S_0$. Fix $s$ in the range $0\leqslant s<S_1$ and consider the number 
\[
\gamma:= F^{(S)} (s\alpha).
\]
The analytic argument (Schwarz's Lemma) provides an upper bound for $|\gamma|$, using the fact that $F$ has at least $SS_1$ zeroes. Liouville's estimate then implies $\gamma=0$.

\medskip
\noindent
{\bf Step 4: \em conclusion}
\\
All derivatives of $F$ at $z=0$ (and also at $z=s\alpha$, $0\leqslant s<S_1$) are $0$, hence the function $F$ is the zero function, which is the desired contradiction with the algebraic independence of $z$ and $\rme^z$.

 \section{Gel'fond's solution of Hilbert's 7th problem}\label{S:G}

The solution of Hilbert's 7th problem has been given independently by Gel'fond \cite{Gelfond} and Schneider \cite{Schneider1}:

\begin{theorem}[Gel'fond--Schneider]\label{Theorem:GS}
Let $\alpha$ and $\beta$ be two complex numbers. Assume $\alpha\not=0$ and $\beta\not\in\Q$. Let $\log\alpha$ be a nonzero logarithm of $\alpha$, meaning $\rme^{\log\alpha}=\alpha$. Define $\alpha^\beta:=\rme^{\beta\log\alpha}$. Then one at least of the three numbers $\alpha$, $\beta$, $\alpha^\beta$ is transcendental.
\end{theorem}

This result implies the transcendence of $2^{\sqrt 2}$, of $\rme^\pi$ (take for instance $\alpha=1$, $\log\alpha=2\rmi \pi$, $\beta=-\rmi/2$), of $\log 2/\log 3$ (take $\alpha=3$, $\beta=\log 2/\log 3$).

The sketch of proof by Gel'fond of Theorem \ref{Theorem:GS} is similar to the one in \S~\ref{S:HL}. One replaces the two functions $z$, $\rme^z$ with $\rme^z$, $\rme^{\beta z}$, which are algebraically independent since $\beta$ is irrational. Also one replaces the points $s\alpha$ of \S~\ref{S:HL} with $s\log \alpha$.

Hence the auxiliary function $F$ is now a linear combination of the functions 
\[
f_{\tau,t}(z)=(\rme^z)^\tau(\rme^{\beta z})^t =\rme^{(\tau +t\beta ) z},
\]
say
\[
F(z)=\sum_{\tau=0}^{T_0-1} \sum_{t=0}^{T_1-1}a_{\tau,t}f_{\tau,t}(z).
\]
 Assume that the three numbers $\alpha$, $\beta$ and $\alpha^\beta$ are algebraic, with $\beta$ irrational. The number field $K$ is now $\Q(\alpha, \beta,\alpha^\beta)$.

Again we select parameters $S_0, S_1,T_0, T_1$. Assuming
 \[
 T_0T_1\geqslant 2[K:\Q]S_0S_1,
 \]
 we deduce from the Thue--Siegel Lemma that there exists a nonzero polynomial $P\in\Z[X,Y]$ such that the function $F(z):=P(\rme^{ z},\rme^{\beta z})$ satisfies 
\begin{equation}\label{Equation:G}
F^{(\sigma)}(s\log \alpha)=0
\end{equation}
for $0\leqslant \sigma<S_0$ and $0\leqslant s<S_1$. By induction, combining Schwarz's Lemma with Liouville's inequality, one proves that equations \eqref{Equation:G} hold for $0\leqslant s<S_1$ and all $\sigma\geqslant 0$. The contradiction follows.

 References are \cite[Ch.~3]{LN402}, \cite[\S~2.1]{Queens}.

 \section{Schneider--Lang Theorem}\label{S:SL}

An entire function $f$ is said to be of finite order if there exists $\varrho>0$ and $R_0>0$ such that $|f|_R\leqslant \rme^{R^{\varrho}}$ for all $R>R_0$. In other terms 
\[
 \limsup_{R\to\infty} \frac {\log\log |F|_R}{\log R}<\infty.
\]

The arguments of sections \ref{S:HL} and \ref{S:G} yield the following statement

 \begin{theorem}[Schneider--Lang]\label{Theorem:SchneiderLang}
 Let $f_1,f_2$ be two entire functions of finite order which are algebraically independent over $\Q$. Let $K$ be a number field. Assume that the derivatives $f'_1$ and $f'_2$ belong to $K[f_1,f_2]$. Then the set
\[
 \{w\in\C\; \mid \; f_1(w) \text{ and } f_2(w) \text{ belong to } K\}
 \]
 is finite.
\end{theorem}
 
The more general version of this Theorem in 
\cite[Ch. III Th.~1]{Lang}, \cite[\S~3.3]{LN402}, \cite[Th.~2.2.1]{Queens}, \cite[Ch.~IV]{GL326}, 
where $f_1$ and $f_2$ are allowed to be meromorphic functions of finite order, has many consequences, in particular in connection with algebraic groups \cite[Ch.~III]{Lang}, \cite[Ch.~5]{Asterisque}, \cite[\S~2.3]{BW}. Also there are generalizations to functions of several variables. For cartesian products, they are due to Schneider and Lang -- see \cite[Ch. IV Th.~1]{Lang}, \cite[Ch.~4]{GL326}. A deeper result involving hypersurfaces, suggested by M. Nagata, has been proved by E. Bombieri; see 
\cite[\S~10.2]{Queens}, \cite[Ch.~5]{Asterisque}.

\medskip
 \noindent
 {\bf Sketch of proof of Theorem \ref{Theorem:SchneiderLang}.}

Let $\varrho_1$, $\varrho_2$ and $C$ satisfy
\[
 \log|f_i|_R\leqslant C R^{\varrho_i} 
\] 
for $i=1,2$ and for all sufficiently large $R$. Assume $w_1,\dots,w_m$ satisfy
\[
 f_i(w_j)\in K \text{ for } i=1,2 \text{ and } j=1,\dots,m.
 \]
 We will prove

\begin{equation}\label{Equation:boundform}
 m\leqslant (\varrho_1+\varrho_2)\delta
\end{equation}
with $\delta:=[K:\Q]$.

 Let $S$ be a sufficiently large integer. We introduce the following parameters $T_1$ and $T_2$:
 \[
 T_1=S^{\frac {\varrho_2} {\varrho_1+\varrho_2} }\sqrt{\log S},
 \quad
 T_2=S^{\frac {\varrho_1} {\varrho_1+\varrho_2} }\sqrt{\log S}.
\]
We need to select these parameters at the beginning of the proof, but this choice is dictated by constraints which occur during the proof. 

The first step is to prove the existence of a nonzero polynomial in $\Z[X_1,X_2]$, of degree $<T_1$ in $X_1$ and $< T_2$ in $X_2$, such that the function 
$F:=P(f_1,f_2)$ satisfies the $mS$ equations 
\begin{equation}\label{Equation:SL}
F^{(\sigma)}(w_j)=0
\end{equation}
for $0\leqslant \sigma<S$ and $1\leqslant j\leqslant m$. Since 
\[
T_1T_2=S\log S,
\]
the existence of $P$ is a consequence of linear algebra as soon as $\log S>\delta m$. For sufficiently large $S$, Thue--Siegel Lemma implies the existence of $P$ with an upper bound for its height:
\[
\log \rmH(P) =O(S).
\]
Since the functions $f_1$ and $f_2$ are algebraically independent, the numbers \eqref{Equation:SL} for $1\leqslant j\leqslant m $ and $\sigma\geqslant 0$ are not all zero. Let $U$ be the least integer such that the function $F$ satisfies \eqref{Equation:SL} for $0\leqslant \sigma<U$ and $1\leqslant j\leqslant m$. From the first step we infer $U\geqslant S$. Since $U$ is minimal, there exists $j_0$ which satisfies $1\leqslant j_0\leqslant m$ such that the number 
\[
\gamma:=F^{(U)}(w_{j_0})
\]
is not $0$. 
Using the inductive assumption together with Schwarz's Lemma with $R=U^{\frac 1 {\varrho_1+\varrho_2} }$, we obtain 
\[
\log |\gamma|\leqslant \left(1-\frac m {\varrho_1+\varrho_2} +o(1)\right)U\log U.
\]
Liouville's estimate for the nonzero element $\gamma$ of $K$ yields 
\[
\log |\gamma|\geqslant -(\delta-1+o(1))U\log U.
\]
The conclusion \eqref{Equation:boundform} follows. 

This is a sketch of the proof in \cite[Ch.3]{LN402} where all details are given. 
See also \cite[Ch. III Th.~1]{Lang} and \cite[Ch. 6]{Baker}.

\section{Schneider's method}\label{S:S}
The proof by Schneider of Theorem \ref{Theorem:GS} is slightly different from Gel'fond's proof in \S~\ref{S:G} as it does not use differential equations.
We argue by contradiction: let again $K$ be the number field $\Q(\alpha,\beta,\alpha^\beta)$. Notice that the two functions $z$ and $\alpha^z$ are algebraically independent since $\log\alpha\not=0$. The auxiliary function $G$ is now a linear combination of the functions 
\[
g_{\sigma,s}(z)=z^{\sigma}\alpha^{sz},
\]
(where, of course, $\alpha^z=\rme^{z\log \alpha}$),
say
\[
G(z)=\sum_{\sigma=0}^{S_0-1} \sum_{s=0}^{S_1-1}b_{\sigma,s}g_{\sigma,s}(z)
\]
with $b_{\sigma,s}\in\Z$. The first step is the existence of a nonzero polynomial 
\[
P(X,Y)=\sum_{\sigma=0}^{S_0-1} \sum_{s=0}^{S_1-1}b_{\sigma,s} X^\sigma Y^s \in\Z[X,Y]
\]
such that the auxiliary function $G(z)=P(z,\alpha^z)$ satisfies
\begin{equation}\label{Equation:SchneiderMethod}
G(\tau+t\beta)=0
\end{equation}
for $0\leqslant \tau<T_0$, $0\leqslant t<T_1$. This follows from the Thue--Siegel Lemma as soon as we request $S_0,S_1,T_0,T_1$ to satisfy 
\[
S_0S_1\geqslant 2[K:\Q]T_0T_1.
\]
By induction on $T'_0+T'_1$, one proves that equations \eqref{Equation:SchneiderMethod} are satisfied for $0\leqslant \tau<T'_0$, $0\leqslant t<T'_1$. For proving that a number of the form $\gamma=G(\tau+t\beta)$ is zero, one uses the analytic estimate given by Schwarz Lemma to show that $|\gamma|$ is small, thanks to the inductive assumptions. One concludes the inductive process by means of Liouville's inequality. Therefore equations \eqref{Equation:SchneiderMethod} are satisfied for all $\tau\geqslant 0$, and $t\geqslant 0$.

The last step is to deduce that $G=0$. There are several ways of achieving this goal, the easiest is to remark that in the inductive argument in the second step, we show that $|G|_r$ is small for arbitrarily large values of $r$, hence $G=0$. 

Two references are \cite[Ch.~2]{LN402} and \cite[\S~3.1]{Queens}. 

\section{Further results}\label{S:FurtherResults}

The reader is invited to use the previous ideas and solve the following exercise:

\begin{quote}
\em Let $x_1,\dots,x_d$ be $\Q$--linearly independent real numbers and $y_1,\dots,y_\ell$ be also $\Q$--linearly independent real numbers. Assume $\ell d>\ell+d$. Then one at least of the $\ell d$ numbers $\rme^{x_iy_j}$ is not a rational number.
\end{quote}

The proof does not use any Liouville's inequality other than the trivial one that a nonzero rational integer has absolute value at least $1$. Using some Liouville's inequality, one reaches the conclusion that one at least of these numbers is transcendental (the $x_i$ and $y_j$ may be complex numbers): this is the so--called {\em Six Exponential Theorem} \cite[Ch. II \S~1]{Lang}, \cite[\S~1.3]{GL326}.

The {\em Four Exponentials Conjecture} is to reach the same conclusion under the weaker assumption $\ell d \geqslant \ell+d$, which means in the case $\ell=d=2$. As pointed out in \cite[Remark, Ch. II \S~1]{Lang}, {\em the same pattern of proof does not give the desired concluslon (it just misses)}. 

We observe that the methods of Gel'fond and Schneider follow the same sketch of proof. There is something more that they share in common: the numbers which occur in Gel'fond's proof and in Schneider's proof are the same. Indeed, with 
\[
f_{\tau,t}(z)=\rme^{(\tau +t\beta ) z} \text{ and } 
g_{\sigma,s}(z)=z^{\sigma}\alpha^{sz}
\] 
as in \S~\ref{S:G} and \S~\ref{S:S}, we have 
\[
f^{(\sigma)} _{\tau,t}(s\log \alpha)=
(\tau+t\beta)^\sigma \alpha^{\tau s}(\alpha^\beta)^{st}
=g_{\sigma,s}(\tau+t\beta).
\] 
This is the source of the duality described in \cite{GL326}. 

In transcendence proofs given in this reference, auxiliary polynomials are replaced with interpolation determinants, an idea which arises in the work by M.~Laurent \cite{Laurent}. 
Consider for instance the above proof in \S~\ref{S:HL} of Hermite--Lindemann Theorem \ref{Theorem:HL}. The algebraic independence of the functions $z$ and $\rme^z$ implies that for a suitable choice of the values for $(\sigma,s)$ and $(\tau,t)$ with $\sigma\geqslant 0$, $s\geqslant 0$, $\tau\geqslant 0$, $t\geqslant 0$, some matrix with entries 
\[
\left(\frac{\rmd}{\rmd z} \right)^\tau (z^\sigma \rme^s)(t\alpha)
\]
has maximal rank. The determinant of this matrix is a nonzero algebraic number. An analytic estimate, based on Schwarz Lemma in a single variable, shows that its absolute value is small. Again we get a contradiction with Liouville's estimate. 

With interpolation determinants, the above mentioned duality between the methods of Gel'fond and Schneider amounts to transposing a matrix. 
For proofs using an interpolation determinant, see \cite[\S~6]{Laurent} for the six exponentials Theorem, \cite[\S~2.3]{GL326} for Schneider's method, \cite[\S~2.4]{GL326} for Gel'fond's method, \cite[\S~2.6]{GL326} for the Hermite--Lindemann Theorem. 

The above sketches of proofs in sections \ref{S:HL}, \ref{S:G} and \ref{S:S} are by contradiction. Since we start with hypotheses that we prove at the end not to be compatible, no step of the proof can be used. However, in practice, one does not argue by contradiction, and an important feature of the method is that it yields effective estimates which have a lot of applications. The first estimates obtained by this approach were due to A.O. Gel'fond, who proved non trivial lower bounds for non vanishing quantities of the form $\Lambda:=\beta\log\alpha_1-\log\alpha_2$ when $\alpha_1$, $\alpha_2$ and $\beta$ are algebraic numbers. It is remarkable that the lower bounds reached by transcendence arguments are valid only under the assumption $\Lambda\not=0$; such estimates are specially useful when $\beta$ is rational: transcendence proofs give lower bounds for non vanishing algebraic numbers of the form $|\alpha_1^{b_1}\alpha_2^{b_2}-1|$. 

Gel'fond's estimates are now a special case of Baker's results \cite[Ch.~2]{Baker} which yield nontrivial explicit lower bounds for so--called linear forms in logarithms
\[
\Lambda:=\beta_0+\beta_1\log\alpha_1+\cdots+\beta_n\log\alpha_n
\]
(with $\alpha_i$ and $\beta_j$ algebraic), 
only assuming $\Lambda\not=0$.

Baker's transcendence result is the following
(\cite[Ch. 2]{Baker}, \cite[Ch. 2]{BW}): 

\begin{theorem}[Baker]\label{Theorem:Baker} 
If $\log \alpha_1,\dots,\log\alpha_n$ are logarithms of algebraic numbers which are linearly independent over $\Q$, then the numbers $1,\log \alpha_1,\dots,\log\alpha_n$ are linearly independent over the field of algebraic numbers. 
\end{theorem}

Hermite--Lindemann Theorem \ref{Theorem:HL} is the case $n=1$, while Gel'fond--Schneider Theorem \ref{Theorem:GS} is the linear independence over the field of algebraic numbers of $\log \alpha_1,\log\alpha_2$. 

Baker's proofs rest on a generalisation of Gel'fond's method. Assuming
\[
\beta_0+\beta_1\log\alpha_1+\cdots+\beta_{n-1}\log\alpha_{n-1}=\log\alpha_n,
\]
he works with $n+1$ functions 
\[
z_0,\rme^{z_1},\dots,\rme^{z_{n-1}},\rme^{\beta_0\beta_1z_1+\cdots+\beta_{n-1}z_{n-1}}
\]
of $n$ variables $z_0, z_1,\dots, z_{n-1}$. The auxiliary function is evaluated, together with its derivatives, at the points multiples of $(1,\log\alpha_1,\dots,\log\alpha_{n-1})$. 

References for proofs using auxiliary functions are given in \cite[\S~10.1.1]{GL326} and \cite[\S~4.1]{Queens} for the homogeneous case $\beta_0=0$, \cite[Ch. 2]{Baker}, \cite[\S~10.1.2 and \S~10.3]{GL326} and \cite[Ch. 2]{BW} for the general case.
Proofs using an interpolation determinant are given in 
\cite[Ch.II \S~1]{LILAN} for the homogeneous case, \cite[Ch.XI]{LILAN} and
\cite[\S~10.2]{GL326} for the general case.

Since all points in Baker's proof lie on the complex line $z(1,\log\alpha_1,\dots,\log\alpha_{n-1})$, the same proof may be described using only functions of a single variable
\[
z,\alpha_1^z,\dots,\alpha_n^z;
\]
see \cite[Ch.~8]{LN402}. In Baker's method, a Schwarz Lemma for functions of a single variable suffices. 

The {\em dual} proof, generalizing Schneider's method to several variables, deals with $n+1$ functions 
\[
z_0,z_1,\dots, z_{n-1},\rme^{z_0}\alpha_1^{z_1}\cdots \alpha_{n-1}^{z_{n-1}}
\]
of $n$ variables $z_0, z_1,\dots, z_{n-1}$. Only the powers of one derivative, $\partial/\partial z_0$, are needed (and no derivative at all in the homogeneous case $\beta_0=0$). The auxiliary function is evaluated at the points in
\[
\Z^n+\Z(\beta_0,\dots,\beta_{n-1}).
\]
With Schneider's method, several variables cannot be avoided, unless $\beta_0=0$ and $n=2$. 
A reference is \cite[\S~6.3]{GL326} for the homogeneous case $\beta_0=0$ and \cite[\S~9.1.1]{GL326} for the general case.

 There is a similar story for the $p$--adic transcendence theory: see for instance \cite[Appendix I]{Asterisque} by D.~Bertrand.


\end{document}